\documentclass[12pt]{article}
\usepackage{amsfonts}
\usepackage{amssymb}
\usepackage{amsmath}
\usepackage{mathrsfs}
\usepackage{CJK}

\begin{document}
\title{ Products of redial derivative and integral-type operators from Zygmund spaces to Bloch spaces}
\author{\large Ning Xu
\footnote{Corresponding author, e-mail: gx899200@126.com.}\\
\small Department of Mathematics and Science, Huai Hai Institute
of Technology,\\
\small Jiangsu, Lianyungang 222005, P. R. China}
\date{}
 \maketitle
 \vskip -1cm
\vskip 3mm

{\bf Abstract.} Let $H(\mathbb{B})$ denote the space of all holomorphic functions on the unit ball
$\mathbb{B}\in \mathbb{C}^n$. In this paper we investigate the boundedness and compactness of the products of radial
derivative operator and the following integral-type operator
$$
I_\varphi^g f(z)=\int_0^1 \Re f(\varphi(tz))g(tz)\frac{dt}{t},\ z\in\mathbb{B}
$$
where $g\in H(\mathbb{B}), g(0)=0$, $\varphi$ is a holomorphic self-map of $\mathbb{B}$,\ between Zygmund spaces and
Bloch spaces.

{\bf Keywords:} radial derivative operator; integral-type operator; Zygmund space; Bloch space

\vskip 3mm

{\bf 1. Introduction}

Let $H(\mathbb{B})$ denote the space of all holomorphic functions on the unit ball $\mathbb{B}\subset\mathbb{C}^n$.
Let $z=(z_1,\ldots,z_n)$ and $w=(w_1,\ldots,w_n)$ be points in the complex vector space $\mathbb{C}^n$ and
$<z,w>=z_1\overline{w_1}+\ldots+z_n\overline{w_n}$. Let
$$
\Re f(z)=\sum z_j\frac{\partial f}{\partial z_j}(z)
$$
stand for the radial derivative of $f\in H(\mathbb{B})$[1]. It is easy to see that, if $f\in H(\mathbb{B}),
f(z)=\sum_\alpha a_\alpha z^\alpha,$ where $\alpha$ is a multi-index, then $\Re f(z)=\sum_\alpha|\alpha|a_\alpha z^\alpha$.
We write $\Re^mf =\Re(\Re^{m-1}f)$.

The Bloch space $\mathcal{B}(\mathbb{B})$ is the space of all $f\in H(\mathbb{B})$ such that
$$
\|f\|_\mathcal{B}:=|f(0)|+\sup\limits_{z\in\mathbb{B}}(1-|z|^2)|\nabla f(z)|<\infty,
$$

The little Bloch space $\mathcal{B}_0(\mathbb{B})$ is the space of all $f\in H(\mathbb{B})$ such that
$$
\lim\limits_{|z|\rightarrow 1}(1-|z|^2)|\nabla f(z)|=0.
$$

It is well known that $f\in \mathcal{B}(\mathbb{B})$ if and only if
$$
b(f):=\sup\limits_{z\in\mathbb{B}}(1-|z|^2)|\Re f(z)|<\infty,
$$
and that $f\in \mathcal{B}_0$ if and only if $\lim\limits_{|z|\rightarrow 1}(1-|z|^2)|\Re f(z)|=0$. Moreover,
the following asymptotic relation holds[2]:
$$
\|f\|_\mathcal{B}\asymp |f(0)|+b(f).
$$

Let $\mathcal{Z}$ denote the class of all $f\in H(\mathbb{B})$, such that
$$
\sup\limits_{z\in\mathbb{B}}(1-|z|^2)|\Re^2 f(z)|<\infty. \ \ \ \ \ \ \ \ \ \ \ \ \ \ \ \ \ \ \ \ \ \ \ \ \ \ \ \ \ \ \ \ \ \ \ \ (1)
$$
Therefore, $\mathcal{Z}$ is called the Zygmund class. Under the natural norm
$$
\|f\|_\mathcal{Z}:=|f(0)|+|f'(0)|+\sup\limits_{z\in\mathbb{B}}(1-|z|^2)|\Re^2 f(z)|<\infty. \ \ \ \ (2)
$$
$\mathcal{Z}$ becomes a Banach space. Zygmund class with this norm will be called the Zygmund space.

The little Zygmund space $\mathcal{Z}_0$ denote the closure in $\mathcal{Z}$ of the set of
all polynomials. From Theorem 7.2 of[3], we see that
$$
f\in \mathcal{Z}_0 \Leftrightarrow \lim\limits_{z\in\mathbb{B}}(1-|z|^2)|\Re^2 f(z)|=0.\ \ \ \ \ \ \ \ \ \ \ \ \ \ \ \ \ \ \ \ \ (3)
$$

Suppose that $g\in H(\mathbb{B}), g(0)=0$ and $\varphi$ is a holomorphic self-map of $\mathbb{B}$,
then an integral-type operator, denote by $I_\varphi^g$ on $H(\mathbb{B})$, is defined as follows:
$$
I_\varphi^g f(z)=\int_0^1 \Re f(\varphi(tz))g(tz)\frac{dt}{t},\ \ g\in H(\mathbb{B}),\ z\in\mathbb{B}  \ \ \ \ \ \ \ \ \ (4)
$$
Operator (4) is related to the following operators
$$
T_g(f)(z)=\int_0^1 f(tz)\Re g(tz)\frac{dt}{t},\ \ I_g(f)=\int_0^1 \Re f(tz)g(tz)\frac{dt}{t}.  \ \ \ \ \ \ \ (5)
$$
acting on $H(\mathbb{B})$, introduced in [4] and studied in [5-10], as well as the operator $T_g$ introduced in [11]acting on holomorphic functions on the unit polydisc (see,also[12],[13],as well as [14] for a particular case of the operator).
One of motivations for introducing operator $I_\varphi^g$ stems from the operator introduced in [15].
Some characterizations of the boundedness and compactness of these and some other integral-type operators mostly in
$\mathbb{C}^n$, can be found, for example, in[4,6,7-9,15-31].

In this paper we study the boundedness and compactness of products of $\Re$ and $I_\varphi^g$ between Zygmund space and
Bloch spaces on the unit Ball.

Throughout this paper, constants are denoted by $C$,they are positive and may differ from one occurence to the other.\\
{\bf 2. Auxiliary results}

{\bf Lemma 1$^{[19]}$.} Let $\Re$ be the radial derivative operator. The product of $\Re$ and $I_\varphi^g$
$$
\Re[ I_\varphi^g(f)](z)=\Re f(\varphi(z))g(z)\ \ \ \ \ \ \ \ \ \ \ \ \ \ \ \ \ \ \ \ \ (6)
$$

{\bf Lemma 2$^{[17]}$.} Suppose $f\in\mathcal{Z}$.The following statements are true.\\
(a). There is a positive constant $C$ independent of $f$ such that
$$
|\Re f(z)|\leq C\|f\|_\mathcal{Z}\ln\frac{e}{1-|z|^2}.\ \ \ \ \ \ \ \ \ \ \ \ \ \ \ \ \ \ \ \ \ (7)
$$
(b). There is a positive constant $C$ independent of $f$ such that
$$
\|f\|_\infty =\sup_{|z|<1}|f(z)|\leq C\|f\|_\mathcal{Z}.\ \ \ \ \ \ \ \ \ \ \ \ \ \ \ \ \ \ \ \ \ (8)
$$

For studying the compactness of the operator $\Re I_\varphi^g:\mathcal{Z}\rightarrow\mathcal{B},$
we need the following Lemma. The proof of the lemma is standard, hence we omit the details.

{\bf Lemma 3.} Assume that $g\in H(\mathbb{B})$, $\varphi$ be a holomorphic self-map of $\mathbb{B}$.
Let $T=\Re I_\varphi^g$ ,then $T:\mathcal{Z}(or\mathcal{Z}_0 )\rightarrow\mathcal{B} $ is compact if and only if $T$
is bounded and for any bounded sequence $(f_k)_{k\in\mathbb{N}}$ in $\mathcal{Z}(or\mathcal{Z}_0 )$
 which converges to 0 uniformly on compact subsets of $\mathbb{B}$, $Tf_k\rightarrow 0$ as $k\rightarrow\infty$.

{\bf Lemma 4$^{[17]}$.} A closed set K in $\mathcal{B}_0$ is compact if and only if it is bounded and satisfies
$$
\lim\limits_{|z|\rightarrow1}\sup\limits_{f\in K}(1-|z|^2)|\Re f(z)|=0.\ \ \ \ \ \ \ \ \ \ \ \ \ \ \ \ \ \ \ \ \ (9)
$$
{\bf 3.The boundedness and compactness of $\Re I_\varphi^g:\mathcal{Z}(\mathcal{Z}_0)
\rightarrow\mathcal{B}(\mathcal{B}_0)$}\\
{\bf Theorem 1.} Let $\varphi$ be a holomorphic self-map of $\mathbb{B}$.
Then the following statements are equivalent.

(a) $\Re I_\varphi^g:\mathcal{Z}\rightarrow\mathcal{B}$ is bounded;

(b) $\Re I_\varphi^g:\mathcal{Z}_0 \rightarrow$ $\mathcal{B}$ is bounded;

(c) $$
\sup\limits_{z\in\mathbb{B}}\frac{(1-|z|^2)|\Re\varphi(z)||g(z)|}{1-|\varphi(z)|^2}<\infty,\ \ \ \ \ \ \ \ \ \ \ \ \ \ \ \ \ \ \ \ \ \ \ \ \ (10)
$$
and
$$
\sup\limits_{z\in\mathbb{B}}(1-|z|^2)|\Re g(z)|\ln\frac{e}{1-|\varphi(z)|^2}<\infty.\ \ \ \ \ \ \ \ \ \ \ \ \ \ \ (11)
$$
{\bf Proof.}(a)$\Rightarrow$(b)  This implication is obvious.

(b)$\Rightarrow$(c) Assume that $\Re I_\varphi^g$:$\mathcal{Z}_0$ $\rightarrow$ $\mathcal{B}$ is boundedness,
i.e.,there exists a constant $C$ such that
$$
\|\Re I_\varphi^g(f)\|_\mathcal{B}\leq C\|f\|_\mathcal{Z},
$$
for all $f\in\mathcal{Z}_0$. Taking the functions $f_j(z)=z_j\in\mathcal{Z}_0$ and
$f_j(z)=z_j-z_j^2\in\mathcal{Z}_0,j=1,2,\cdots,n$,we get
$$
\sup\limits_{z\in\mathbb{B}}(1-|z|^2)|\varphi_j(z)||\Re \varphi(z) g(z)+\Re g(z)|<\infty,\ \ \ \  (12)
$$
and
$$
\sup\limits_{z\in\mathbb{B}}(1-|z|^2)|(\varphi_j(z)-4\varphi^2_j(z))
(\Re\varphi(z) g(z)+\Re g(z))+2\varphi^2_j(z)\Re g(z)|<\infty.\ (13)
$$
Using (12) and the boundedness of functions $\varphi_j$, we have that
$$
\sup\limits_{z\in\mathbb{B}}(1-|z|^2)|\Re \varphi(z) g(z)+\Re g(z)|<\infty.\ \ \ \  (14)
$$
Then with (13),(14) and the boundedness of functions $\varphi_j$, we have that

$$
\sup\limits_{z\in\mathbb{B}}(1-|z|^2)|\Re g(z)|<\infty ,
\ \ \sup\limits_{z\in\mathbb{B}}(1-|z|^2)|\Re\varphi(z) g(z)|<\infty  \ \ \ \ (15)
$$

Set
$$
h(\zeta)=(\zeta-1)[(1+\ln\frac{1}{1-\zeta})^2+1],  \zeta\in \mathbb{C},
$$
and
$$
h_a(z)=\frac{h(<z,a>)}{|a|^2}(\ln\frac{1}{1-|a|^2})^{-1},
$$
for $a\in\mathbb{B}\backslash \{0\}$. It is known that $h_a(z)\in\mathcal{Z}_0$(see[17]).Since
$$
\Re h_a(z)=\frac{<z,a>}{|a|^2}(\ln\frac{1}{1-<z,a>})^2(\ln\frac{1}{1-|a|^2})^{-1},
$$
and
$$
\Re^2h_a(z)=\Re h_a(z)+\frac{2<z,a>^2}{|a|^2(1-<z,a>)}(\ln\frac{1}{1-<z,a>})(\ln\frac{1}{1-|a|^2})^{-1},
$$
for$|\varphi(z)|>\sqrt{1-1/e}$ we have
\begin {eqnarray*}
C\|\Re I_\varphi^g\|_\geq||\Re I_\varphi^g(h_{\varphi(z)})\|_\mathcal{B}
&\geq& (1-|z|^2)\ln\frac{1}{1-|\varphi(z)|^2}\Re g(z)|\\
&-&\frac{2(1-|z|^2)}{1-|\varphi(z)|^2}|\varphi(z)|^2||\Re\varphi(z)||g(z)|\\
&-&(1-|z|^2)\ln\frac{1}{1-|\varphi(z)|^2}|\Re\varphi(z)||g(z)|.
\end {eqnarray*}
Hence
\begin {eqnarray*}
(1-|z|^2)\ln\frac{1}{1-|\varphi(z)|^2}|\Re g(z)|&\leq & C+\frac{2(1-|z|^2)}{1-|\varphi(z)|^2}|\varphi(z)|^2||\Re\varphi(z)||g(z)|\\
&+&(1-|z|^2)\ln\frac{1}{1-|\varphi(z)|^2}|\Re\varphi(z)||g(z)|\\
&\leq &C+(2+e)\frac{1-|z|^2}{1-|\varphi(z)|^2}|\Re\varphi(z)||g(z)|,\ \ \ \ \ (16)
\end {eqnarray*}
which using the fact of $(1-|\varphi(z)|^2)\ln\frac{1}{1-|\varphi(z)|^2}\leq e$.\\
For $|a|>\sqrt{1-1/e}$, set
$$
f_a(z)=\frac{h(<z,a>)}{|a|^2}(\ln\frac{1}{1-|a|^2})^{-1}-\int_0^1\ln\frac{1}{1-<tz,a>}\frac{dt}{t}.
$$
Then $f_a\in \mathcal{Z}_0$.It is easy to see that
$$
\Re f_a(z)=\frac{<z,a>}{|z|^2}(\ln\frac{1}{1-<z,a>})^2(\ln\frac{1}{1-|a|^2})^{-1}-\ln\frac{1}{1-<z,a>},
$$
\begin {eqnarray*}
\Re^2f_a(z)&=&\Re f_a(z)+\frac{2<z,a>^2}{|a|^2(1-<z,a>)}(\ln\frac{1}{1-<z,a>})(\ln\frac{1}{1-|a|^2})^{-1}\\
&-&\frac{<z,a>}{1-<z,a>}+\ln\frac{1}{1-<z,a>}.
\end {eqnarray*}
Therefore
\begin {eqnarray*}
C\|\Re I_\varphi^g\|\geq\|\Re I_\varphi^g(f_{\varphi(z)})\|_\mathcal{B}
&=&\sup\limits_{z\in \mathbb{B}}(1-|z|^2)|\Re^2(I_\varphi^gf_{\varphi(z)})(z)|\\
&=&(1-|z|^2)(\frac{|\varphi(z)|^2}{1-|\varphi(z)|^2}+\ln\frac{1}{1-|\varphi(z)|^2})|\Re \varphi(z)||g(z)|\\
&\geq&(1-|z|^2)(\frac{|\varphi(z)|^2}{1-|\varphi(z)|^2}+1)|\Re \varphi(z)||g(z)|\\
&=&\frac{1-|z|^2}{1-|\varphi(z)|^2}|\Re \varphi(z)||g(z)|.\ \ \ \ \ \ \ \ \ \ \ \ \ \ \ \ (17)
\end {eqnarray*}
On the other hand, from (15) we have that
$$
\sup\limits_{|\varphi(z)|\leq\sqrt{1-1/e}}(1-|z|^2)|\Re g(z)|\ln\frac{1}{1-|\varphi(z)|^2}\leq
\sup\limits_{|\varphi(z)|\leq\sqrt{1-1/e}}(1-|z|^2)|\Re g(z)|<\infty.\ \ \ (18)
$$
Hence from (15),(16),(17)and (18),we obtain (11).Further, from (17), we have
$$
\sup\limits_{|\varphi(z)|>\sqrt{1-1/e}}\frac{1-|z|^2}{1-|\varphi(z)|^2}|\Re \varphi(z)||g(z)|\leq C.\ \ \ \ \ \ \ \ \ \ \ \ \ \ \ \ \ \ \ \ \ \ \ \ \ \ (19)
$$
On the other hand, from (15) we have that
$$
\sup\limits_{|\varphi(z)|\leq\sqrt{1-1/e}}\frac{1-|z|^2}{1-|\varphi(z)|^2}|\Re \varphi(z)||g(z)|\leq
e, \  \sup\limits_{|\varphi(z)|\leq\sqrt{1-1/e}}(1-|z|^2)|\Re \varphi(z)||g(z)|<\infty.\ \ \ \ \ \ \ \ \ \ \ (20)
$$
Combining (19)and (20),(10) follows.

{\bf Theorem 2.} Let $\varphi$ be a holomorphic self-map of $\mathbb{B}$. Then the following
statements are equivalent.

(a)$\Re I_\varphi^g:\mathcal{Z}\rightarrow\mathcal{B}$ is compact;

(b)$\Re I_\varphi^g:\mathcal{Z}_0\rightarrow\mathcal{B}$ is compact;

(c)$\Re I_\varphi^g:\mathcal{Z}\rightarrow\mathcal{B}$ is bounded,
$$
\lim_{|\varphi(z)|\rightarrow1}\frac{1-|z|^2}{1-|\varphi(z)|^2}|\Re \varphi(z)||g(z)|=0,\ \ \ \ \ \ \ \ \ \ \ \ \ \ \ \ \ \ \ \ \ \ (21)
$$
and
$$
\lim_{|\varphi(z)|\rightarrow1}(1-|z|^2)|\Re g(z)|\ln\frac{e}{1-|\varphi(z)|^2}=0.\ \ \ \ \ \ \ \ \ \ \ \ \ \ \ \ (22)
$$

{\bf Proof.} (a)$\Rightarrow$ (b)  This is obvious.

 (b)$\Rightarrow$ (c) Assume that $\Re I_\varphi^g:\mathcal{Z}_0\rightarrow\mathcal{B}$ is compact,
then it is clear that $\Re I_\varphi^g:\mathcal{Z}_0\rightarrow\mathcal{B}$ is bounded. By theorem 1,
we know that $\Re I_\varphi^g:\mathcal{Z}\rightarrow\mathcal{B}$ is bounded. Let $(z^k)_{k \in\mathbb{N}}$
be a sequence in $\mathbb{B}$ such that $|\varphi(z^k)|\rightarrow1$ as $k\rightarrow\infty$
 and $\varphi(z^k)\neq 0, k\in\mathbb{N}$. Set
$$
h_k(z)=\frac{h(<z,\varphi(z^k)>)}{|\varphi(z^k)|^2}(\ln\frac{1}{1-|\varphi(z^k)|^2})^{-1},  k\in\mathbb{N}.
$$
Then from the proof of theorem 1,we see that $h_k\in\mathcal{Z}_0$,for each $k\in\mathbb{N}$.
Moreover $h_k\rightarrow 0$ uniformly on compact subsects of $\mathbb{B}$, as $k\rightarrow\infty$.

Since $\Re I_\varphi^g:\mathcal{Z}_0\rightarrow\mathcal{B}$ is compact,by Lemma 3
$$
\lim_{k\rightarrow\infty}\|\Re [I_\varphi^g(h_k)]\|_\mathcal{B}=0.
$$
On the other hand, similar to the proof of Theorem 1, we have
\begin {eqnarray*}
0 \leftarrow\|\Re I_\varphi^g(h_k)\|_\mathcal{B}
&\geq& (1-|z^k|^2)\ln\frac{1}{1-|\varphi(z^k)|^2}|\Re g(z^k)|\\
&-&\frac{2(1-|z^k|^2)}{1-|\varphi(z^k)|^2}|\varphi(z^k)|^2\Re\varphi(z^k)||g(z^k)|\\
&-&(1-|z^k|^2)\ln\frac{1}{1-|\varphi(z^k)|^2}|\Re\varphi(z^k)||g(z^k)|\\
&=&(1-|z^k|^2)\ln\frac{1}{1-|\varphi(z^k)|^2}|\Re g(z^k)|\\
&-&M_1\frac{(1-|z^k|^2)}{1-|\varphi(z^k)|^2}|\Re\varphi(z^k)||g(z^k)|,
\end {eqnarray*}
where $M_1:=2|\varphi(z^k)|^2-(1-|\varphi(z^k)|^2)\ln\frac{1}{1-|\varphi(z^k)|^2}$.

From this we obtain
$$
\lim\limits_{k\rightarrow\infty}(1-|z^k|^2)\ln\frac{1}{1-|\varphi(z^k)|^2}|\Re g(z^k)|=
\lim\limits_{k\rightarrow\infty}\frac{(1-|z^k|^2)}{1-|\varphi(z^k)|^2}|\Re\varphi(z^k)||g(z^k)|=0,\ \ \ \ \ (23)
$$
if one of these two limits exists,which use the case of $$
\lim\limits_{k\rightarrow\infty}[2|\varphi(z^k))^2|+(1-|\varphi(z^k)|^2)\ln\frac{1}{1-|\varphi(z^k)|^2}]=2.
$$

Next, set
\begin {eqnarray*}
f_k(z)&=&\frac{h(<z,\varphi(z^k)>)}{|\varphi(z^k)|^2}(\ln\frac{1}{1-|\varphi(z^k)|^2})^{-1}\\
&-&\int_0^1\ln^3\frac{1}{1-<tz,\varphi(z^k)>}\frac{dt}{t}(\ln\frac{1}{1-|\varphi|(z^k)|^2})^{-2}.
\end {eqnarray*}
Since $\Re I_\varphi^g:\mathcal{Z}_0\rightarrow\mathcal{B}$
is compact,we have
$\|\Re I_\varphi^g(f_k)\|_{\mathcal{B}}\rightarrow0$ as $k\rightarrow\infty.$
Thus,similar to the proof of Theorem 1,when$\sqrt{1-\frac{1}{e}}<|\varphi(z^k)|<1$
\begin {eqnarray*}
0\leftarrow\|\Re I_\varphi^g(f_k)\|_\mathcal{B}
&\geq&(1-|z^k|^2)|\ln\frac{1}{1-|\varphi(z^k)|^2}-\frac{|\varphi(z^k)|^2}{1-|\varphi(z^k)|^2}||\Re \varphi(z^k)||g(z^k)|\\
&\geq&(1-|z^k|^2)|\frac{|\varphi(z^k)|^2}{1-|\varphi(z^k)|^2}||\Re \varphi(z^k)||g(z^k)|\\
&-&(1-|z^k|^2)\ln\frac{1}{1-|\varphi(z^k)|^2}|\Re \varphi(z^k)||g(z^k)|\\
&\geq&(1-\frac{1}{e})\frac{1-|z^k|^2}{1-|\varphi(z^k)|^2}|\Re \varphi(z^k)||g(z^k)|\\
&-&\frac{1-|z^k|^2}{1-|\varphi(z^k)|^2}(1-|\varphi(z^k)|^2)\ln\frac{1}{1-|\varphi(z^k)|^2}|\Re \varphi(z^k)||g(z^k)|\\
&=&M_2\frac{1-|z^k|^2}{1-|\varphi(z^k)|^2}|\Re \varphi(z^k)||g(z^k)|,
\end {eqnarray*}
where$ M_2:=1-\frac{1}{e}-(1-|\varphi(z^k)|^2)\ln\frac{1}{1-|\varphi(z^k)|^2}$.

Hence
$$
\lim\limits_{k\rightarrow\infty}\frac{1-|z^k|^2}{1-|\varphi(z^k)|^2}|\Re \varphi(z^k)||g(z^k)|=
\lim\limits_{k\rightarrow\infty}(1-|z^k|^2)\ln\frac{1}{1-|\varphi(z^k)|^2}|\Re g(z^k)|=0.\ (24)
$$

From (24) easily following that $\lim\limits_{k\rightarrow\infty}(1-|z^k|^2)|\Re g(z^k)|=0$,which altogether imply (21)and (22).

(c)$\Rightarrow$ (a)
$$
C_1=(1-|z|^2)|\Re \varphi(z)||g(z)|<\infty,\ \ \ C_2=(1-|z|^2)|\Re g(z)|<\infty.\ \ \ \ \ \ \ \ \ \ \ (25)
$$

For every $\varepsilon>0$,there is a $\delta\in (0,1)$, such that
$$
\frac{1-|z|^2}{1-|\varphi(z)|^2}|\Re \varphi(z)||g(z)|<\varepsilon,\ \ \
(1-|z|^2)|\Re g(z)|\ln\frac{e}{1-|z|^2}<\varepsilon.\ \ \ \ \ \ \ \ \ \ \ \ \ (26)
$$

Assume that $(f_k)_{k\in\mathbb{N}}$ is a sequence in $\mathcal{Z}$
such that $\sup\limits_{k\in\mathbb{N}}\|f_k\|_{\mathcal{Z}}\leq L$ and $f_k$ converges to 0
uniformly on compact subsets of $\mathbb{B}$ as $k\rightarrow\infty$.
Let $K=\{z\in\mathbb{B}:|\varphi(z)|\leq\delta\}$. Then by (25) and (26) ,we have that
\begin {eqnarray*}
\sup\limits_{z\in\mathbb{B}}(1-|z|^2)|(\Re^2(I_{\varphi}^g(f_k))(z)|
&=&\sup\limits_{z\in\mathbb{B}}(1-|z|^2)|\Re^2f_k(\varphi(z))\Re\varphi(z)g(z)+\Re f_k(\varphi(z))\Re g(z)|\\
&\leq&\sup\limits_{z\in\mathbb{B}}(1-|z|^2)|\Re^2f_k(\varphi(z))\Re\varphi(z)g(z)|\\
&+&\sup\limits_{z\in\mathbb{B}}(1-|z|^2)|\Re f_k(\varphi(z))\Re g(z)|\\
&\leq&\sup\limits_{z\in K}(1-|z|^2)|\Re^2f_k(\varphi(z))\Re\varphi(z)g(z)|\\
&+&\sup\limits_{z\in K}(1-|z|^2)|\Re f_k(\varphi(z))\Re g(z)|\\
&+&\sup\limits_{z\in\mathbb{B}\backslash K}(1-|z|^2)|\Re^2f_k(\varphi(z))\Re\varphi(z)g(z)|\\
&+&\sup\limits_{z\in\mathbb{B}\backslash K}(1-|z|^2)|\Re f_k(\varphi(z))\Re g(z)|\\
&\leq&\sup\limits_{z\in K}(1-|z|^2)|\Re^2f_k(\varphi(z))\Re\varphi(z)g(z)|\\
&+&\sup\limits_{z\in K}(1-|z|^2)|\Re f_k(\varphi(z))\Re g(z)|\\
&+&\sup\limits_{z\in\mathbb{B}\backslash K}\frac{1-|z|^2}{1-|\varphi(z)|^2}|\Re\varphi(z)g(z)|\|f_k\|_{\mathcal{Z}}\\
&+&C \sup\limits_{z\in\mathbb{B}\backslash K}\ln\frac{e}{1-|\varphi(z)|^2}|\Re g(z)|\|f_k\|_{\mathcal{Z}}\\
&\leq&C_1\sup\limits_{z\in K}|\Re^2f_k(\varphi(z))|+C_2\sup\limits_{z\in K}|\Re f_k(\varphi(z))|+(C+1)\varepsilon\|f_k\|_{\mathcal{Z}}\\
\end {eqnarray*}

Hence
\begin {eqnarray*}
\|\Re I_{\varphi}^g(f_k)\|_\mathcal{B}&\leq& C_1\sup\limits_{z\in K}|\Re^2f_k(\varphi(z))|+C_2\sup\limits_{z\in K}|\Re f_k(\varphi(z))|\\
&+&(C+1)\varepsilon\|f_k\|_{\mathcal{Z}} +|f'_k(\varphi(0))||\varphi'(0)|
\end {eqnarray*}

Since $(f_k)_{\mathbb{N}}$ converges to 0 uniformly on compact subsets of $\mathbb{B}$ as $k\rightarrow\infty$, Cauchy's estimate gives that
$\Re f_k\rightarrow 0$ and $\Re^2 f_k\rightarrow 0$ as $k\rightarrow\infty$ on compact subsets of $\mathbb{B}$.
Hence, letting $k\rightarrow\infty$, we obtain
$$
\lim\limits_{k\rightarrow\infty}\|\Re I_{\varphi}^g(f_k)\|_\mathcal{B}=0.
$$

{\bf Theorem 4.} Let $\varphi$ be a holomorphic self-map of $\mathbb{B}$. Then $\Re I_\varphi^g:
\mathcal{Z}_0\rightarrow\mathcal{B}_0$ is bounded if and only if $\Re I_\varphi^g:\mathcal{Z}_0\rightarrow\mathcal{B}$ is bounded
$$
\lim_{|z|\rightarrow1}(1-|z|^2)|\Re g(z)|=0,\ \ \ \ \ \ \ \ \ \ \ \ \ \ \ \ \ \ \ (27)
$$
$$
\lim_{|z|\rightarrow1}(1-|z|^2)|\Re \varphi(z)||g(z)|=0.\ \ \ \ \ \ \ \ \ \ \ \ \ (28)
$$
{\bf Proof:} Assume that $\Re I_\varphi^g:\mathcal{Z}_0\rightarrow\mathcal{B}_0$ is bounded. Then, it is clear that
$\Re I_\varphi^g:\mathcal{Z}_0\rightarrow\mathcal{B}$ is bounded.Taking the function $f_j(z)=z_j$ and $f_j(z)=z_j-z_j^2,j=1,2,\cdots,n,$
we obtain (27),(28).

Assume that $\Re I_\varphi^g:\mathcal{Z}_0\rightarrow\mathcal{B}$ is bounded and (27),(28) holds. Then for each polynomial $p$,
we have that
\begin {eqnarray*}
(1-|z|^2)|\Re^2( I_\varphi^gp)(z))|&\leq &(1-|z|^2)|\Re^2p(\varphi(z))||\Re \varphi(z)||g(z)|\\
&+ &(1-|z|^2)|\Re p(\varphi(z))||\Re g(z)|,\ \ \ \ \ \ \ \ \ \ \ \ \ (29)
\end {eqnarray*}
from (27),(28) it follows that $\Re I_\varphi^g p \in\mathcal{B}_0$. Since the set of all polynomials is dense in
$\mathcal{Z}_0$, we have that for every $f\in\mathcal{Z}_0$, there is a sequence of polynomials $(p_n)_{n\in\mathbb{N}}$
 such that $\|f-p_n\|_\mathcal{Z}\rightarrow0$ as $n\rightarrow\infty$.
 Hence
 $$
 \|\Re I_{\varphi}^g(f)-\Re I_{\varphi}^g(p_n)\|_\mathcal{B}
 \leq \|\Re I_{\varphi}^g\|_{\mathcal{Z}_0\rightarrow\mathcal{B}}\|f-p_n\|_\mathcal{Z}\rightarrow 0
 $$
as $n\rightarrow\infty$. Since the operator $\Re I_{\varphi}^g:\mathcal{Z}_0\rightarrow\mathcal{B}$ is bounded,
hence $\Re I_{\varphi}^g(\mathcal{Z}_0)\subseteq \mathcal{B}_0$.

{\bf Theorem 5.} Let $\varphi$ be a holomorphic self-map of $\mathbb{B}$.Then the following statements are
equivalent.

(a)$\Re I_\varphi^g:\mathcal{Z}\rightarrow\mathcal{B}$ is compact;

(b)$\Re I_\varphi^g:\mathcal{Z}_0\rightarrow\mathcal{B}_0$ is compact;

(c)
$$
\lim_{|z|\rightarrow1}\frac{1-|z|^2}{1-|\varphi(z)|^2}|\Re \varphi(z)||g(z)|=0,\ \ \ \ \ \ \ \ \ \ \ \ \ (30)
$$
and
$$
\lim_{|z|\rightarrow1}(1-|z|^2)|\Re g(z)|\ln\frac{e}{1-|\varphi(z)|^2}=0.\ \ \ \ \ \ \ \ \ \ \ \ \ (31)
$$

{\bf Proof:} (a)$\Rightarrow$(b). It is clear.

(b)$\Rightarrow$(c).Assume that $\Re I_\varphi^g:\mathcal{Z}_0\rightarrow\mathcal{B}_0$ is compact,then
$\Re I_\varphi^g:\mathcal{Z}_0\rightarrow\mathcal{B}_0$ is bounded.From the proof of Theorem 4,we known that
$$
\lim_{|z|\rightarrow1}(1-|z|^2)|\Re g(z)|=0,
$$
$$
\lim_{|z|\rightarrow1}(1-|z|^2)|\Re \varphi(z)||g(z)|=0,
$$

Hence, if $\|\varphi\|<1$,
$$
\lim_{|z|\rightarrow1}\frac{1-|z|^2}{1-|\varphi(z)|^2}|\Re \varphi(z)||g(z)|\leq
\frac{1}{1-\|\varphi\|_\infty}\lim_{|z|\rightarrow1}(1-|z|^2)|\Re g(z)|=0,
$$
$$
\lim_{|z|\rightarrow1}(1-|z|^2)|\Re g(z)|\ln\frac{e}{1-|\varphi(z)|^2}\leq
\ln\frac{e}{1-\|\varphi\|^2_\infty}\lim_{|z|\rightarrow1}(1-|z|^2)|\Re g(z)|=0.
$$
from which the result follows in this case.

Assume $\|\varphi\|=1$. Let $(\varphi(z^k))_{k\in\mathbb{N}}$ be a sequence such that $|\varphi(z^k)|\rightarrow 1$
as $k\rightarrow\infty$. Since $\Re I_\varphi^g:\mathcal{Z}_0\rightarrow\mathcal{B}$ is compact, by Theorem 2,
$$
\lim_{|\varphi(z)|\rightarrow1}\frac{1-|z|^2}{1-|\varphi(z)|^2}|\Re \varphi(z)||g(z)|=0,\ \ \ \ \ \ \ \ \ \ \ \ \ (32)
$$
and
$$
\lim_{|\varphi(z)|\rightarrow1}(1-|z|^2)|\Re g(z)|\ln\frac{e}{1-|\varphi(z)|^2}=0.\ \ \ \ \ \ \ \ \ \ \ \ \ (33)
$$
It is not difficult to see that (28),(32) implies(30). Similar, (27) and (33) imply (31).

(c)$\Rightarrow$(a). Let $f\in \mathcal{Z}$, we have
$$
(1-|z|^2)|\Re^2(I_{\varphi}^g(f))(z)|\leq(\frac{1-|z|^2}{1-|\varphi(z)|^2}|\Re \varphi(z)||g(z)|+
(1-|z|^2)\ln\frac{e}{1-|\varphi(z)|^2}|\Re g(z)|\|f\|_{\mathcal{Z}}.
$$
Taking the supremum in this inequality over all $f\in \mathcal{Z}$ such that $\|f\|_{\mathcal{Z}}\leq 1$.
Letting $|z|\rightarrow1$ and using (30),(31)
$$
\lim\limits_{\|z\|\rightarrow1}\sup\limits_{\|f\|_{\mathcal{Z}}\leq 1}(1-|z|^2)|\Re^2(I_{\varphi}^g(f))(z)|=0.
$$
Using Lemma 3,we obtain that the operator $\Re I_\varphi^g:\mathcal{Z}\rightarrow\mathcal{B}_0$ is compact.

 \vspace{0.2cm}
\begin{center}{REFERENCES}
\end{center}

\vspace{-0.3cm}
\begin{enumerate}

\item  W.Rudin, Function Theory in the Unit Ball of $\mathbb{C}^n$, Spring-Verlay, New York, 1980.

\item  D.Clahane, S.Stevi$\acute{c}$, Norm equivalence and composition operators between
Bloch/Lipschitz spaces of the unit ball, J.Inequal. Appl. 2006(2006).Article ID 61068, 11 pp.

\item  K.Zhu, Spaces of Holomorphic Functions in the Unit Ball, Springer, New York, 2005.

\item  Z.Hu, Extended Ces$\grave{a}$ro operators on the Bloch spaces in the unit ball of $\mathbb{C}^n$,
Acta. Math. Sci. Ser. B Engl. Ed.23(4)(2003)561-566.

\item  S.Stevi$\acute{c}$, On an integral operator on the unit ball in $\mathbb{C}^n$,J.Inequeal.Appl.\\2005(1)81-88.

\item  N.Xu, Extended Ces$\grave{a}$ro operators on $\mu$-Bloch spaces in $\mathbb{C}^n$,
J. of Math. Research  and  Exposition, 29(5)(2009)913-922.

\item Z.Hu, Extended Ces$\grave{a}$ro operators on mixed norm spaces, Proc. Amer. Math. Soc. 131(7)(2003)2171-2179.

\item  Z.Hu, Extended Ces$\grave{a}$ro operators on  Bergman spaces, J. Math. Anal. Appl. 296(2004)435-454.

\item  S. Li, S.Stevi$\acute{c}$,  Ces$\grave{a}$ro-type operators on some spaces of analytic functions on the unit ball,Appl. Math. and Comput. 208(2009)378-388.

\item  S.Stevi$\acute{c}$, Ces$\grave{a}$ro averaging operators,Math.Nachr.248-249(2003)185-189.

\item  S.Stevi$\acute{c}$, Boundedness and compactness of an integral operator on a weighted  space on the polydisc,Indian J. Pure Appl. Math. 37(6)(2006) 343-355.

\item  S.Stevi$\acute{c}$, Boundedness and compactness of an integral operator on mixed norm spaces on polydisc, Sibirsk. Mat. Zh. 48(3)(2007)694-706.

\item  D.C. Chang, S.Stevi$\acute{c}$, The generalized Ces$\grave{a}$ro operator on the unit polydisk,
Taiwan. J. Math. 7(2)(2003)293-308.

\item D.C. Chang, S.Stevi$\acute{c}$, Estimates of an integral operator on function spaces,
Taiwan. J. Math. 7(3)(2003)423-432.

\item  S. Li, S.Stevi$\acute{c}$, Generalized composition operators on Zygmund spaces and Bloch type spaces, J. Math. Anal. Appl. 338(2008)1282-1295.

\item  D.C. Chang, S.Stevi$\acute{c}$, On some integral operators on the unit polydisk and the unit ball,
Taiwan. J. Math. 11(5)(2007)1251-1286.

\item  S. Li, S.Stevi$\acute{c}$, Riemann-Stieltjes operators between mixed norm spaces, Indian J. Math. 50(1)(2008)177-188.

\item  S. Li, S.Stevi$\acute{c}$, Products of composition and differentiation operators from Zygmund spaces to Bloch spaces and Bers spaces,Appl. Math. and Comput. 217(2010)3144-3154.

\item  S.Stevi$\acute{c}$, On an integral operator between Bloch-type spaces on the unit ball, Bull. Sci. Math. 134(2010)329-339.

\item  S.Stevi$\acute{c}$, S.I.Ueki, Integral-type operators acting between weighted-type spaces on unit ball,Appl.Math.Comput. 215(2009)2464-2471.

\item  K.Avetisyan, S.Stevi$\acute{c}$, Extended Ces$\grave{a}$ro operator between different Hardy spaces, Appl.Math.Comput.(2009)346-350.

\item  W.Yang, On an integral-type operator between Bloch-type spaces,\\Apll.Math.Comput.215(3)(2009)954-960.

\item  X.Zhang, Weighted composition operators between $\mu$-Bloch spaces on the unit ball,Sci.China(Ser.A)48(2005)1349-1368.

\item  X.Zhu, Integral-type operators from iterated logarithmic Bloch spaces to Zygmund-type spaces,Appl.Math.Comput.215(2009)1170-1175.

\item  S.Stevi$\acute{c}$, Generalized composition operators from logarithmic Bloch spaces to mixed-norm spaces,Util.Math.77(2008)167-172.

\item  S.Stevi$\acute{c}$, On an integral operator from the Zygmund space to the Bloch-type space on the unit ball, Glasg.J.Math.51(2009)275-287.

\item  S.Stevi$\acute{c}$, On an integral operator from the logarithmic Bloch-type and  mixed-norm spaces to Bloch-type spaces, Nonlinear Anal. TMA 71(2009)6323-6342.

\item  S.Stevi$\acute{c}$, Weighted composition operators from the logarithmic weighted-type space to the weighted Bergman space in $\mathbb{C}^n$, Appl.Math.Comput.\\216(2010)924-928.

\item  S.Stevi$\acute{c}$, On a new operator from H$^\infty$ to the Bloch-type space on the unit ball,Util.Math.77(2008)257-263.

\item  S.Stevi$\acute{c}$, On a new operator from the logarithmic Bloch space to the Bloch-type space on the unit ball,Appl.Math.Comput.206(2008)313-320.

\item  S.Stevi$\acute{c}$, Products of integral-type operators and composition operators from the mixed-norm space to Bloch-type spaces,Siberian J.Math.50(4)(2009)726-736.

\end{enumerate}

\end{document}